\documentclass[10pt, xcolor = {dvipsnames}
]{amsart}
\usepackage[dvipsnames]{xcolor}
\usepackage{tikz} 
\usepackage{pgfplots}
\usepackage{amsmath, amsthm, amssymb,mathrsfs,amsfonts,}
\usepackage{mathtools}
\usepackage{relsize}
\usepackage{enumitem}
\usepackage[utf8]{inputenc}
\usepackage{pgffor}
\usepackage{color}
\usepackage{epstopdf}
\usepackage{graphicx}

\usepackage{varwidth}
\usepackage{verbatim}

\usepackage[foot]{amsaddr}

\usepackage{mathtools}      
\mathtoolsset{showonlyrefs} 

\pgfplotsset{compat=1.17}

\usepackage[
]
{hyperref}
\hypersetup{
	colorlinks=false,       
	linkcolor=blue,          
	citecolor=red,        
	filecolor=magenta,      
	urlcolor=cyan           
}

\newtheorem{thm}{Theorem}[section]

\newtheorem{lem}[thm]{Lemma}

\theoremstyle{definition}

\theoremstyle{remark}
\newtheorem{rk}[thm]{Remark}

\numberwithin{equation}{section}

\definecolor{ao}{rgb}{0, 0.5, 0}

\begin{document}

\title[Bourgain's counterexample for the sequential convergence]{Bourgain's counterexample in the sequential convergence problem for the Schr\"odinger equation}

\author{Chu-Hee Cho}
\address{Department of Mathematical Sciences and RIM, Seoul National University, Republic of Korea
		  e-mail: {\tt akilus@snu.ac.kr}
		}

\author{Daniel Eceizabarrena}
\address{BCAM - Basque Center for Applied Mathematics, Bilbao, Spain, and Department of Mathematics, 
University of Massachusetts Amherst, USA
		  e-mail: {\tt deceizabarrena@bcamath.org}
		}

\thanks{}
\date{\today}

\begin{abstract}
We study the 
problem of pointwise convegence for the Schr\"odinger operator on $\mathbb R^n$ along time sequences.  
We show that the sharp counterexample to the sequential Schr\"odinger maximal estimate given recently by Li, Wang and Yan based in the construction by Lucà and Rogers can also be achieved with the construction of Bourgain, 
and we extend it to the fractal setting.
\end{abstract}
\maketitle

\section{Introduction}
The Schr\"odinger operator
\begin{equation}\label{Schrodinger operator}
e^{it\Delta}f(x) = \int_{\mathbb{R}^n} \widehat f(\xi) \, e^{ 2\pi i (x \cdot \xi + t|\xi|^2) } \,  d\xi
\end{equation}
gives the solution to the free Schr\"odinger equation on $\mathbb R^n$ with initial datum $f$.  
In 1980, Carleson \cite{Carl} asked for the smallest Sobolev regularity $s$ such that 
\begin{align}\label{conv_original}
\lim_{t \rightarrow 0} e^{it\Delta}f = f  \quad \text{ a.e.}, \quad \text{ for all } f \in H^s(\mathbb R^n).
\end{align}
This problem is now solved; the critical Sobolev exponent is $n/(2(n+1))$ \cite{Carl,DK,B2,DGL,DZ}.
Interesting variations of this problem have been studied since then.  
In this note we focus on two of them: the fractal and the sequential problems. 

In the fractal problem,  the convergence property is strengthened to $\mathcal H^{\alpha}$-a.e., where  $\alpha \leq n$ is given and $\mathcal H^{\alpha}$ is the $\alpha$-Hausdorff measure. 
That is, we look for 
\begin{equation}
s_c(\alpha) = \inf \{ \, s\ge 0 \,  :  \,  \lim_{t \rightarrow 0} e^{it\Delta}f = f \quad \mathcal H^{\alpha} \text{- a.e.}, \, \, \,  \forall  f\in H^s(\mathbb R^n) \, \}.  
\end{equation}
\v Zubrini\'c \cite{Z} and Barcel\'o et al. \cite{BBCR} showed $s_c(\alpha) = (n-\alpha)/2$ when $0<\alpha \le n/2$,
but the problem is open when $n/2 < \alpha \le n$.
The best results are the upper bound  
\begin{equation}\label{DuZhang}
s_c(\alpha) \leq \frac n{2(n+1)} + \frac n{2(n+1)} (n-\alpha), \quad \text{ for } \quad  \frac{n + 1}2 \leq \alpha \le n 
\end{equation}
by Du--Zhang \cite{DZ}, and a sawtooth-like lower bound\footnote{We refer to \cite[Theorem 1.1]{EP_V}  for the explicit expression of this exponent.} $s_{\text{E-PV}}$ by  Eceizabarrena--Ponce-Vanegas \cite[Theorem 1.1]{EP_V}  that recently improved the previous best lower bound by Luc\`a--Rogers \cite{LR} and Luc\`a--Ponce-Vanegas \cite{LP-V}:
\begin{equation} \label{Best_Lower_Bound}
\frac n{2(n+1)} + \frac {n-1}{2(n+1)} (n-\alpha)  \leq s_{\text{E-PV}}(\alpha) \leq  s_c(\alpha).
\end{equation}

In the sequential problem, 
given a sequence of times $t_k \to 0$ in a class $T$ that measures its speed of convergence, 
we look for the smallest $s$ such that 
\begin{align}\label{conv_sequential}
\lim_{k \to \infty} e^{it_k \Delta} f = f  \quad \text{ a.e.  \quad for all } (t_k) \in T,   \qquad \text{ for all } f \in H^s(\mathbb R^n).
\end{align}
In \cite{DS}, Dimou and Seeger characterized this regularity for $n = 1$ in terms of Lorentz spaces
$T = \ell^{r,\infty}$, 
which for $0   < r  < \infty$ are defined by 
\begin{equation}
(t_k) \in \ell^{r,\infty}(\mathbb N) 
\qquad \Longleftrightarrow \qquad 
 \sup_{\delta>0} \delta^r \#\{ k \in \mathbb N: |t_k| \geq \delta \} < \infty. 
\end{equation}
Denoting the corresponding critical exponent 
by 
\begin{equation}
 s_c(\ell^{r,\infty})
=  \inf \{ \, s\ge 0 \,  :  \,  \lim_{k \to \infty} e^{it_k\Delta}f = f \quad \text{a.e.}, \, \forall (t_k) \in \ell^{r,\infty}, \, \,  \forall  f\in H^s(\mathbb R^n)\}, 
\end{equation}
we now know from \cite{DS,CKKL,LWY2} that
\begin{equation}\label{Critical_Sequential_Exponent}
s_c(\ell^{r,\infty}) = \min \Big(\frac{n}{2(n+1)}, \frac{rn}{r(n+1)+n} \Big),
\end{equation}
which shows a non-trivial change of behavior at $r = n/(n+1)$. 


In this paper we work with the combination of these two problems
and the corresponding fractal sequential exponent
\begin{equation}\label{Critical_Exponent_Fractal_Sequential}
 s_c(\alpha, \ell^{r,\infty})
=  \inf \{ \, s\ge 0 \,  :  \,  \lim_{k \to \infty} e^{it_k\Delta}f = f \quad \mathcal H^{\alpha}\text{-a.e.}, \, \forall (t_k) \in \ell^{r,\infty}, \, \,  \forall  f\in H^s(\mathbb R^n)\}. 
\end{equation}
By standard arguments,
upper bounds for \eqref{Critical_Exponent_Fractal_Sequential} follow from maximal estimates 
\begin{equation}\label{fractal_maximal_estimate}
\int_{B(0,1)} \sup_{k}|e^{it_k\Delta}f(x) |^2 d\mu(x) \lesssim C_\mu \|f\|_{H^s}^2, 
\quad \forall (t_k) \in \ell^{r,\infty}, \quad \forall f \in H^s, \quad \forall \mu \in \mathcal M_\alpha,  
\end{equation}
where $\mathcal M_\alpha$ is the set of $\alpha$-dimensional measures, that is, 
non-negative Borel measures for which there is a constant $C_\mu$ such that
\begin{align}\label{Def_Measure}
\mu (B(x,r)) \le C_\mu r^\alpha, \quad \forall x \in \mathbb R^n, \quad r > 0.
\end{align}
Let us denote the critical Sobolev exponent for the fractal sequential maximal estimate by
\begin{equation}
s^\text{max}_c(\alpha, \ell^{r, \infty}) 
=   \inf \{ \, s\ge 0 \,  :  \,  \eqref{fractal_maximal_estimate} \text{ holds } \},
\end{equation}
and by $s^\text{max}_c(\alpha)$ we denote the exponent for the non-sequential maximal estimate, where the supremum in \eqref{fractal_maximal_estimate} is taken over $t \in (0,1)$ instead of over a sequence $t_k$. 


In \cite{CKKL},
Cho, Ko, Koh and Lee proved 
\begin{equation}\label{CKKL}
s^\text{max}_c(\alpha) \leq \frac{n - \alpha}{2} + s^* 
\quad \Longrightarrow \quad 
s^\text{max}_c(\alpha, \ell^{r, \infty}) \leq \frac{n - \alpha}{2} + \min \Big(s^* , \frac{2r s^*}{r + 2s^*} \Big). 
\end{equation}
With the upper bounds known for $s^\text{max}_c(\alpha)$, 
one gets the following positive result: 
\begin{thm}\label{thmPositive}
Let $n\geq 1$, $0 \leq \alpha \leq n$ and $0<r<\infty$. Then, 
\begin{itemize}
	\item If $\alpha \leq n/2$, then 
		\begin{equation}
			s^\text{max}_c(\alpha, \ell^{r, \infty}) = \frac{n - \alpha}{2}. 
		\end{equation}
		
	\item If $n/2 < \alpha \leq (n+1)/2$, then 
		\begin{equation}
			s^\text{max}_c(\alpha, \ell^{r, \infty}) 
			 \leq \frac{n - \alpha}{2} + \min \Big(\frac{2\alpha - n}{4} , \frac{r(2\alpha - n)}{2r + 2\alpha - n} \Big).
		\end{equation}
	
	\item If $\alpha \geq (n+1)/2$, then 
	\begin{equation}
		s^\text{max}_c(\alpha, \ell^{r, \infty}) 
		\leq \frac{n - \alpha}{2} + \min \Big(\frac{\alpha}{2(n+1)} , \frac{r\alpha}{r(n+1) + \alpha} \Big).
	\end{equation} 
\end{itemize}
\end{thm}

We briefly explain in Section~\ref{sec:PositiveResult} how to prove this combining \eqref{CKKL} 
with the results by Zubrini\'c \cite{Z},  Barceló et al \cite{BBCR}, and Du and Zhang \cite{DZ}. 
Observe that when  $\alpha \geq (n+1)/2$, 
\begin{equation}
	r < \frac{\alpha}{n+1} \qquad \Longrightarrow \qquad  s^\text{max}_c(\alpha, \ell^{r, \infty})  \leq  \frac{n-\alpha}2 + \frac{r\alpha }{(n+1)r+\alpha},
	\end{equation}
	improving the upper bound by Du-Zhang. However, there is no improvement if
	\begin{equation}
	r \geq \frac{\alpha}{n+1} \qquad \Longrightarrow \qquad  s^\text{max}_c(\alpha, \ell^{r, \infty})   \leq  \frac{n}{2(n+1)}(n+1-\alpha). 
	\end{equation}	
 When $n/2 \leq \alpha \leq (n+1)/2$, we improve the known $n/4$ when $r < (2\alpha - n)/2$.


Regarding lower bounds, when $\alpha = n$ we know the optimal \eqref{Critical_Sequential_Exponent}. 
In this note
we give the first lower bounds for $s_c^\text{max}(\alpha,\ell^{r,\infty})$ when $\alpha  < n$. 

In $n=1$, the typical single wave-packet argument gives:
\begin{thm}\label{main_theorem_1D}
Let $n=1$ and $1/2 \leq \alpha \leq 1$. Then, 
\begin{equation}
s^\text{max}_c(\alpha, \ell^{r, \infty}) \geq \frac{1-\alpha}{2} + \min\Big( \frac{2\alpha - 1}{4}, \frac{r(2\alpha - 1)}{2r + \alpha} \Big).
\end{equation}
\end{thm}

\begin{rk}
Together with Theorem~\ref{thmPositive}, this theorem implies: 
\begin{enumerate}
	\item If $1/2 \leq \alpha \leq 2r$, then $s^\text{max}_c(\alpha, \ell^{r, \infty}) = 1/4$. 
	
	\item If $2r \leq \alpha \leq r + 1/2$, then 
	\begin{equation}
	\frac{1-\alpha}{2} + \frac{r(2\alpha - 1)}{2r + \alpha} \leq s^\text{max}_c(\alpha, \ell^{r, \infty}) \leq \frac14.
\end{equation}		
	
	\item If $r + 1/2 \leq \alpha \leq 1$, then 
	\begin{equation}
	\frac{1-\alpha}{2} + \frac{r(2\alpha - 1)}{2r+\alpha} 
	\leq s^\text{max}_c(\alpha, \ell^{r, \infty}) 
	\leq \frac{1-\alpha}{2} + \frac{r(2\alpha - 1)}{2r+2\alpha - 1}.
	\end{equation}
	
\end{enumerate}
In particular, the critical exponent is  1/4 if $r \geq 1/2$. 
\end{rk}


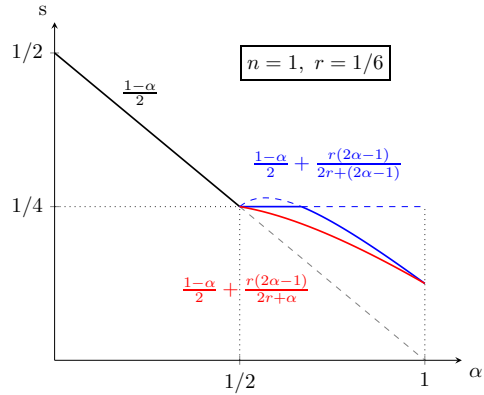
\begin{figure}[h]
\begin{center}
\begin{tikzpicture}[scale=0.79]
\begin{axis}[axis x line=center, axis y line=center, 
xmin=0, xmax=1.1, ymin=0, ymax=0.55, xlabel={$\alpha$}, ylabel={s},  xtick={0, 0.5,1}, xticklabels={$0$, $1/2$, $1$}, ytick={0.25,0.5}, yticklabels={$1/4$,$1/2$},
xlabel style={below right},
ylabel style={above left}
]
\addplot[black, domain=0:0.5, thick] {(1-x)/2};
\addplot[gray, domain=0.5:1, dashed] {(1-x)/2};
\addplot[blue, domain=0.5:2/3, thick] {1/4};
\addplot[blue, domain=2/3:1, dashed] {1/4};
\addplot[black, dotted] coordinates {(0,0.25) (0.5,0.25)};
\addplot[black, dotted] coordinates {(0.5,0) (0.5,0.25)};
\addplot[black, dotted] coordinates {(1,0) (1,0.25)};
\addplot[blue,  domain=0.5:2/3, dashed] {(1-x)/2+(2*x-1)/(2 + 6*(2*x-1))};
\addplot[blue,  domain=2/3:1, thick] {(1-x)/2+(2*x-1)/(2 + 6*(2*x-1))};
\addplot[red, domain=0.5:1, thick] {(1-x)/2+(2*x-1)/(6*x+2)};
\node[above=4.5cm, right=1cm] {$\frac{1-\alpha}{2}$};
\node[above=3.3cm, right=3.2cm, blue] {$\frac{1-\alpha}{2}+\frac{r(2\alpha-1)}{2r+(2\alpha-1)}$};
\node[above=1.2cm, right=2cm, red] {$\frac{1-\alpha}{2}+\frac{r(2\alpha-1)}{2r+\alpha}$};
\node[above=5cm, right=3cm, black] {$\boxed{n=1, \, \, r=1/6}$};
\end{axis}
\end{tikzpicture}
\caption{Upper (blue) and lower (red) bounds when $n=1$.}
\end{center}
\end{figure}

%

When $n \geq 2$, to get \eqref{Critical_Sequential_Exponent}
Li, Wang and Yan \cite{LWY2} adapted the ergodic counterexample of Lucà and Rogers \cite{LR}. 
We show here that the same result can also be obtained via Bourgain's argument in \cite{B2}, and we generalize it to the fractal setting.  
\begin{thm}\label{main theorem}
Let $n\ge 2$. 
\begin{enumerate}[font=\normalfont]
	\item
	\label{MainTheorem_Item1}	
	If $n-1 \leq \alpha \leq n$, 
	\begin{equation}\label{Exponent_Big_Alpha}
	s^\text{max}_c(\alpha, \ell^{r, \infty}) \geq \frac{n-\alpha}2 + \min \Big( \frac{2\alpha - n}{2(n+1)},  \frac{r(2\alpha-n)}{r(n+1)+2\alpha-n} \Big).
	\end{equation}
	\item
	\label{MainTheorem_Item2}	
	 If $n/2 < \alpha \leq n - 1$, 
	\begin{equation}\label{Exponent_Small_Alpha}
	s^\text{max}_c(\alpha, \ell^{r, \infty}) \geq \frac{n-\alpha}2 +  \min \Big( \frac{2\alpha - n}{2(n+1)}, \frac{r\alpha}{r(n+1) + n} \Big). 
	\end{equation}
	
\end{enumerate}
\end{thm}

\begin{rk}
\begin{enumerate}
	\item The exponent in \eqref{Exponent_Big_Alpha} is larger than the one in \eqref{Exponent_Small_Alpha}
	due to a better scaling argument in the range $n-1 \leq \alpha \leq n$. 
	
	\item When $r \geq (2\alpha - n)/(n+1)$, in both cases the exponent is the same as the one in \eqref{Best_Lower_Bound} by Lucà-Rogers and Lucà-Ponce-Vanegas.   
	The exponents are new when $r < (2\alpha - n)/(n+1)$.
	
	

	\item The exponent $s_{\text{E-PV}}$ in \eqref{Best_Lower_Bound} was obtained in \cite{EP_V} by generalizing Bourgain's counterexample with an \textit{intermediate space trick} from \cite{DKWZ}.
	Combining our strategy in this note with this intermediate space trick should give a better lower bound,
	but it will require heavy computations. 
	We prefer to keep this note short to highlight the ideas of our construction, so we will not pursue this further here. 
	
\end{enumerate}
\end{rk}

\begin{figure}[h]
\begin{center}
\begin{tikzpicture}[scale=1]
\begin{axis}[axis x line=center, axis y line=center,
xmin=0, xmax=3.2, ymin=0, ymax=1.6, xlabel={$\alpha$}, ylabel={s}, xtick={0,1.5,2.65,3}, xticklabels={$0$, $n/2$, $n-1$, $n$}, ytick={0.75,1.5}, yticklabels={$n/4$,$n/2$},  
xlabel style={below right},
ylabel style={above left}
]
\addplot[black, domain=0:1.5, thick] {(3-x)/2};
\addplot[gray, domain=1.5:3, dashed] {(3-x)/2};
\addplot[gray, domain=0:1.5, dashed] {3/4};
\addplot[gray, dashed] coordinates {(1.5,0) (1.5,33/40)};
\addplot[gray, dashed] coordinates {(2.65,0) (2.65,0.335)};
\addplot[gray, dashed] coordinates {(3,0) (3,3/20)};
\addplot[red,  domain=1.5:3, thick] {(3-x)/2+(6*x-9)/(32*x-36)};
\draw[color=blue, thick] (1.5,0.79)--(3,3/20);
\node[above=4.5cm, right=1.2cm] {$\frac{n-\alpha}{2}$};
\node[above=1.5cm, right=1.2cm, blue] {$\frac{n-\alpha}{2}+\frac{r\alpha}{r(n+1)+n}$};
\node[above=3cm, right=3.5cm, red] {$\frac{n-\alpha}{2}+\frac{r(2\alpha-n)}{r(n+1)+(2\alpha-n)}$}; 
\end{axis}
\end{tikzpicture}
\caption{Comparison of the two lower bounds in Theorem \ref{main theorem}.} 
\end{center}
\end{figure}
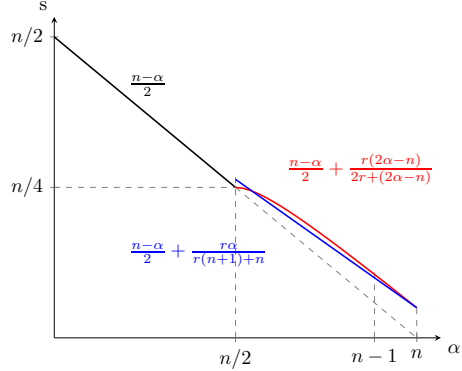

\subsubsection*{Notation}
$A \lesssim B$ denotes $A \leq CB$ for a constant $C>0$. 
We write $A \simeq B$ when $A \lesssim B$ and $B \lesssim A$.
For $\mu \in \mathcal M_\alpha$, we denote the infimum of $C_\mu$ in \eqref{Def_Measure} by $\langle\mu \rangle_\alpha$. 

\section{Brief explanation of Theorem \ref{thmPositive}}\label{sec:PositiveResult}

Theorem~\ref{thmPositive} directly follows from the transference result \eqref{CKKL} proved by Cho et al., which we reproduce here for convenience of the reader:
\begin{equation}\label{CKKL_2}
s^\text{max}_c(\alpha) \leq \frac{n - \alpha}{2} + s^* 
\quad \Longrightarrow \quad 
s^\text{max}_c(\alpha, \ell^{r, \infty}) \leq \frac{n - \alpha}{2} + \min \Big(s^* , \frac{2r s^*}{r + 2s^*} \Big). 
\end{equation}
If $\alpha \leq n/2$, by Barcel\'o et al. \cite{BBCR} \eqref{CKKL_2} holds with $s^*=0$. 
Together with the result by \v Zubrini\'c \cite{Z}, this directly implies that 
\begin{equation}
s_c(\alpha, \ell^{r, \infty}) = s^\text{max}_c(\alpha, \ell^{r, \infty}) = \frac{n - \alpha}{2}, \qquad \forall r >0.
\end{equation}
When $n/2 \leq \alpha \leq (n+1)/2$, the bound by Barceló et al. \cite{BBCR} and the fact that $s_c^\text{max}(\alpha)$ is decreasing imply that $s^\text{max}_c(\alpha) \leq n/4$. 
Hence, \eqref{CKKL_2} holds with $s^* = (2\alpha - n)/4$, so
\begin{equation}
\begin{split}
s^\text{max}_c(\alpha, \ell^{r, \infty}) 
& \leq \frac{n - \alpha}{2} + \min \Big(\frac{2\alpha - n}{4} , \frac{r(2\alpha - n)}{2r + 2\alpha - n} \Big) \\
& = \min \Big(\frac{n}{4} , \frac{n - \alpha}{2} + \frac{r(2\alpha - n)}{2r + 2\alpha - n}  \Big),
\end{split}
\end{equation}
which gives an improvement if $r \leq (2\alpha - n)/2$. 
Finally, by Du-Zhang \cite{DZ}, when $\alpha \geq (n+1)/2$ we get \eqref{CKKL_2} with $s^* = \alpha / (2(n+1))$, 
which implies 
\begin{equation}
s^\text{max}_c(\alpha, \ell^{r, \infty}) 
\leq \frac{n - \alpha}{2} + \min \Big(\frac{\alpha}{2(n+1)} , \frac{r\alpha}{r(n+1) + \alpha} \Big).
\end{equation} 
This improves the Du-Zhang upper bound if $r \leq \alpha /(n+1)$.

\section{Proof of theorem \ref{main theorem} - Dimensions $n \geq 2$}\label{sec:Proof_Main_Thm}


Let $R \gg 1$ be the main parameter. 
We will construct a sequence of functions $f_R$ and measures $\mu_R$ to disprove the maximal estimate \eqref{fractal_maximal_estimate}. 
Let $S,D$ be auxiliary parameters with $1 \leq S \leq R$ and $1 \leq D \leq R$, 
both of which are powers of $R$ that we determine later.  

Let $\phi$ be a smooth function\footnote{With a small abuse of notation, $\phi$ and $\psi$ may represent functions both in $\mathbb R$ and $\mathbb R^{n-1}$.} with $\operatorname{supp} \widehat{\phi} \subset B_1(0)$ and $\phi(0) = 1$.
Let $\psi$ be a positive smooth function with support in 
$\{ 
 1/2 \leq |x| \leq 1 \}$. 
Define $f_R$ by 
\begin{equation}
\widehat{f}_R (\xi) = \widehat{g}_R (\xi_1) \,  \widehat{h}_R (\xi'),
\qquad \xi=(\xi_1, \xi') \in \mathbb R \times \mathbb R^{n-1}, 
\end{equation}
where 
\begin{equation}\label{Def_Of_ID}
\widehat{g}_R(\xi_1)= \frac{1}{S} \,  \widehat{\phi}\Big(\frac{\xi_1-R}{S}\Big),  
 \qquad \qquad
\widehat{h}_R(\xi')= \sum_{m'\in \mathbb Z^{n-1}} \psi\Big( \frac{D}{R} m' \Big) \, \widehat{\phi}(\xi'-Dm').
\end{equation}
It follows that 
$\|g_R\|_2 \simeq S^{-1/2}$ and  $\|h_R\|_2 \simeq (R/D)^{\frac{n-1}2}$, so
\begin{equation}\label{l2norm}
\|f_R\|_2 = \|g_R\|_2\|h_R\|_2 \simeq \frac{1}{S^{1/2}} \, \Big( \frac{R}{D} \Big)^{\frac{n-1}2}.
\end{equation}
To study the evolution of $f_R$, we do so for $g_R$ and $h_R$ separately. 
For $x_1 \in [0,1]$, 
\begin{equation}
\begin{split}
e^{it\Delta}g_R(x_1) &= \int   \frac{1}{S} \,  \widehat{\phi}\Big(\frac{\xi_1-R}{S}\Big) \, e^{2\pi i ( x_1\xi_1 + t\xi_1^2)}  d\xi_1\\
                                 &= e^{2\pi i (x_1R + tR^2)} \, \int  \widehat{\phi}(\eta_1)  \, e^{2\pi i ( S\eta_1(x_1 + 2Rt) + tS^2\eta_1^2) } \,d\eta_1,
\end{split}
\end{equation}
where we changed variables $\xi_1 = R+S\eta_1$. 
For $c>0$ be small enough,
\begin{equation}\label{con1}
 S|x_1+2Rt| \le c, \quad S^2t \le c
 \quad  \Longrightarrow \quad 
 |e^{it\Delta}g_R(x_1)| \simeq 1.
\end{equation} 
On the other hand, 
for $x' \in \mathbb [0,1]^{n-1}$ we write
\begin{equation}\label{evolution2}   
\begin{split}
&
e^{it\Delta}h_R(x') \\
& \qquad  =  \sum_{m'\in \mathbb Z^{n-1}} \psi\Big( \frac{D}{R} m' \Big)  \, \int \widehat{\phi}(\xi'-Dm') \, e^{2\pi i (x'\cdot\xi' + t|\xi'|^2)}\, d\xi' \\
& \qquad = \sum_{m'\in \mathbb Z^{n-1}} \psi\Big( \frac{D}{R} m' \Big) \,  e^{2\pi i (Dx'\cdot m' + tD^2|m'|^2)} \,  \int  \widehat{\phi}(\eta')  e^{2\pi i (\eta' \cdot(x'+2tDm') + t|\eta'|^2)} \, d\eta'.
\end{split}
\end{equation}
Since $|\eta'|\le 1$ and $D|m'| \simeq R$,  we ask
\begin{equation}\label{con2}
t \leq \frac{c}{R}
\end{equation} 
to make the phase inside the integral small. 
Simultaneously, we set
\begin{equation}\label{con3}
D^2t =   \frac {p_1}q,
\qquad 
Dx_i = \frac {p_i}q + \epsilon, \quad  i = 2, \ldots, n
\end{equation}  
with $q \in \mathbb Z$, $p_1, \cdots, p_n \in \mathbb Z$, $\operatorname{gcd}(p_1,q) = 1$ and $|\epsilon| \leq cD/R$. 
Conditions \eqref{con2} and \eqref{con3} allow us to estimate the evolution \eqref{evolution2}
by a Gauss sum, 
that is, 
\begin{align}
|e^{it\Delta}h_R(x')| 
\label{Gauss1}    & \simeq 
\Big| \sum_{m'\in \mathbb Z^{n-1}} \psi\Big( \frac{D}{R} m' \Big) \,  e^{2\pi i \frac{p' \cdot m' + p_1 |m'|^2}{q} } \, e^{2\pi i \epsilon \cdot m'}\Big| \\
\label{Gauss2}    & \simeq \prod_{i=2}^n  \Big( \frac{R}{Dq} \Big) \Big| \sum_{m_i=1}^q  e^{2\pi i \frac{p_i \,  m_i + p_1 m_i^2}{q} } \Big| 
\end{align}
where $p' = (p_2,\cdots, p_n)$ and, slightly abusing notation again, $\epsilon = (\epsilon, \ldots, \epsilon) \in \mathbb R^{n-1}$. 
By properly choosing $p$ and $q$ (for example, choosing $q$ odd or, alternatively, $q$ even and $p_i \equiv q/2 \pmod{2}$), 
each Gauss sum has size $\sqrt{q}$, so 
\begin{equation}\label{other variables}  
|e^{it\Delta}h_R(x')| 
 \simeq  \Big( \frac{R}{Dq} \, \sqrt{q} \Big)^{n-1}  = \Big( \frac{R}{D\sqrt{q}} \Big)^{n-1}. 
\end{equation}
Altogether\footnote{
For a detailed proof of estimates \eqref{Gauss1}, \eqref{Gauss2} and \eqref{other variables}, we refer to \cite[Section 3.1]{EP_V2} or \cite{P}.
}, combining \eqref{l2norm}, \eqref{con1} and \eqref{other variables} we get
\begin{equation}\label{solution_lb}
\frac{|e^{it\Delta}f_R(x)|}{\lVert f_R \rVert_2} 
\simeq S^{1/2} \,  \Big( \frac R{Dq}\Big)^{(n-1)/2}.
\end{equation} 

Let us now discuss the conditions we imposed in \eqref{con1}, \eqref{con2}, \eqref{con3}. 
First, rewrite the second condition in \eqref{con3} and $\epsilon \le cD/R$ as
\begin{equation}\label{x'_ball}
x' \in B_{n-1}\Big(\frac1D\frac {p'}q, \frac cR\Big).
\end{equation}
On the other hand, from the first conditions in \eqref{con1} and \eqref{con3} we get
\begin{equation}\label{x1_ball}
x_1 
\in 
 B_{1}\Big(\frac R{D^2}\frac {p_1}q, \frac cS\Big).
\end{equation}
Hence, for a new parameter $Q$ (also a power of $R$ to be determined later), let 
\begin{equation}\label{baseset}
X_R = \bigcup_{q \simeq Q} \bigcup_p B_1 \Big( \frac R{D^2} \frac{p_1}{q}, \frac cS \Big) \times B_{n-1}\Big(\frac1D \frac{p'}{q}, \frac cR \Big)
\subset [0,1]^n,
\end{equation}
where $q \simeq Q$ means $Q/2 \leq q \leq Q$, possibly with additional restrictions that do not alter substantially the number of $q$ considered\footnote{For example, parity conditions, or imposing $q$ to be prime.}.
Regarding $t$, \eqref{con1} and \eqref{con2} imply
\begin{equation}
t \le c \, \min \Big( \frac{1}{S^2}, \frac{1}{R} \Big). 
\end{equation}
This condition determines two regimes, depending on $S \leq R^{1/2}$ or not. 
In view of the condition $S|x_1 + 2Rt| \leq c$ in \eqref{con1}, it affects the set $X_R$: 
\begin{enumerate}
	\item[(i)] When $S \lesssim R^{1/2}$, then $t \lesssim 1/R$ so all $ |x_1| \lesssim 1$ are allowed. 
	This is the situation of the counterexample in Bourgain's article \cite{B2}. 
	
	 \item[(ii)] When $R^{1/2} \ll S \leq R$, then $t \lesssim 1/S^2 \ll 1/R$,  so 
	 \begin{equation}\label{X_R_Small}
	 t \lesssim \frac{1}{S^2} \ll \frac{1}{R} \quad \Longrightarrow \quad |x_1| \leq \frac{R}{S^2} \ll 1 \quad \Longrightarrow \quad  X_R \subset [0,R/S^2] \times [0,1]^{n-1}. 
	\end{equation}	  
	This situation is the analogue to the one in Li--Wang--Yan \cite{LWY2} to give the optimal counterexample for the sequential convergence.
\end{enumerate}
From now on, we focus on case (ii) above and assume that $R^{1/2} \le S \le R$.
In particular,
from $D^2 t = p_1 / q \lesssim D^2/S^2$
and $|x'| \leq c$
we get the restrictions 
\begin{equation}
p_1 \lesssim \frac{qD^2}{S^2}, \qquad p' \lesssim Dq. 
\end{equation}
Gathering it all, we have shown that 
\begin{equation}\label{Pointwise_Lower_bound}
x \in X_R 
\quad \Longrightarrow \quad 
\exists t(x) \in \Big[\frac{1}{D^2 Q}, \frac{1}{S^2} \Big]
\,  \, : \,  \,  
\frac{|e^{it(x)\Delta}f_R(x)|}{\lVert f_R \rVert_2} \simeq S^{1/2} \,  \Big( \frac R{DQ}\Big)^{(n-1)/2}.
\end{equation}

%
%
%

\subsection{Lebesgue measure of $X_R$}\label{sec:Lebesgue_Measure}
Let us choose the parameters $S,D,Q$ in such a way that $X_R$ has the largest possible measure. 
In view of \eqref{X_R_Small}, we aim for $|X_R| \simeq R/S^2$. 
We first argue heuristically:
\begin{itemize}
	\item[\textbf{Heuristic argument 1.}] It is a consequence of Minkowski's theorem that if $a_1, \ldots, a_n \geq 1$ are such that $a_1 + \cdots + a_n = n+1$, then any $x \in [0,1]^n$ can be approximated as 
	\begin{equation}
	\big| x_i - \frac{p_i}{q} \big| \leq \frac{1}{q^{a_i}}, \qquad i=1, \ldots, n. 
	\end{equation}
	If we rescale $X_R$ in \eqref{baseset} to balls centered at rationals and write the radii as $q^{-a_i}$, we get 
	\begin{equation}
	\frac{D^2}{RS} \simeq \frac{1}{Q^{a_1}}, \qquad \frac{D}{R} \simeq \frac{1}{Q^{a_i}}, \quad i = 2, \ldots, n. 
	\end{equation}
	Thus, asking the Minkowski condition $n+1 = a_1 + \ldots + a_n = a_1 + (n-1)a_2$, 
	\begin{equation}
	\frac{1}{Q^{n+1}} = \frac{1}{Q^{a_1} \, Q^{(n-1)a_2}} =  \frac{D^2}{RS} \, \Big( \frac{D}{R} \Big)^{n-1}
	\qquad \Longleftrightarrow \qquad (DQ)^{n+1} \simeq S \, R^n. 
	\end{equation}

	\item[\textbf{Heuristic argument 2.}] A simpler (but more heuristic) argument is to assume that the union \eqref{baseset} of $X_R$ is disjoint, 
	so that
\begin{equation}\label{Volume}
|X_R| \simeq Q\cdot \frac{D^2Q}{S^2}\cdot(DQ)^{n-1}\cdot\frac1S\frac1{R^{n-1}} \sim \frac{(DQ)^{n+1}}{S^3R^{n-1}}.
\end{equation}	
The condition that we impose is then
\begin{equation}\label{measure_con}
\frac{(DQ)^{n+1}}{S^3R^{n-1}} \simeq \frac{R}{S^2} \qquad \Longleftrightarrow \qquad (DQ)^{n+1} \simeq S R^n,
\end{equation}
which is the same as above. 
\end{itemize}
These heuristics can be made rigorous using Lemma~\ref{Lemma_AnChuPierce} in Appendix~\ref{Appendix_LemmaAnChuPierce}.
Define the dilation $T(x_1,x') = ((D^2/R) x_1, Dx')$, so that 
$|T(X_R)| = (D^2/R) D^{n-1} |X_R|$. 
The set $T(X_R)$ has a 1-periodic structure and it is made of 
$(D^2/S^2) D^{n-1}$
disjoint translations of the set 
\begin{equation}\label{UnitCell}
\widetilde X_R = \bigcup_{q \sim Q} \bigcup_{p=1}^q B_1 \Big( \frac{p_1}{q}, c\frac{D^2}{RS} \Big) \times B_{n-1}\Big( \frac{p'}{q}, c\frac{D}{R} \Big).
\end{equation}
This means that $(D^2/S^2) D^{n-1}  |\widetilde X_R| = |T(X_R)| = (D^2/R) D^{n-1} |X_R|$, so
\begin{equation}\label{X_R_in_terms_of_tildeXR}
 |X_R| =  |\widetilde X_R| \, \frac{R}{S^2}. 
\end{equation}
Now, calling $h_1(Q) = D^2/(RS)$ and $h_2(Q) = D/R$, Lemma~\ref{Lemma_AnChuPierce} implies
\begin{equation}
(QD)^{n+1} \simeq S R^n 
\quad \Longrightarrow \quad 
Q^{n+1} h_1(Q) h_2(Q)^{n-1} \simeq 1 
\quad \Longrightarrow \quad
|\widetilde X_R| \simeq 1,
\end{equation}
and therefore 
\begin{equation}\label{Measure_X_R_Lebesgue}
|X_R| \simeq R/S^2.
\end{equation}
\begin{rk}\label{rk:Small_Measure}
In Bourgain's counterexample, in which $S \lesssim R^{1/2}$, one gets $|X_R| \simeq 1$ for all $R$. 
From it, one can explicitly construct a counterexample to the convergence problem by working with the set of dicergence $X = \limsup_R  X_R$ 
with positive measure $|X| \simeq 1$. 
In this case, however, from \eqref{X_R_Small} it is clear that $X = \limsup_R X_R \subset \{0\} \times [0,1]^{n-1}$, 
so in particular it has measure $|X| = 0$.  
Even worse, $\operatorname{dim}_{\mathcal{H}} X \leq n-1$. 
So even if this disproves the maximal estimate, it is more challenging to explicitly construct a counterexample to the convergence problem,  
at least when $n-1 < \alpha \leq n$; we will not pursue this issue in this note. 
\end{rk}

\subsection{Sequence of times in $\ell^{r,\infty}$} 
The times $t$ we selected in \eqref{con3} belong in  
\begin{equation}
T = \bigcup_R T_R, \quad  \text{ where } \quad T_R = \Big\{  \, \frac{1}{D^2} \frac{p_1}{q} \, : \, q \simeq Q, \, \, 1 \leq p_1 \leq \frac{D^2 q}{S^2}\,  \Big\}. 
\end{equation}
Observe first that $T_R \subset [1/(D^2Q), 1/S^2]$. 
To make all these intervals disjoint, 
we select a sequence of dyadic $R_n$ such that $1/S_{n+1}^2 \ll 1/(D_n^2Q_n)$. 
This way, we can treat the set of times $T$ as a sequence that tends to zero.
By choosing $S,D,Q$ properly, we will make it belong in $\ell^{r,\infty}$. 
For that, we ask
\begin{equation}\label{time sequence}
\delta^r \#\{t : \delta < t < 2 \delta \} \lesssim 1, \qquad \text{ for every } \delta>0.
\end{equation}
Suppose $1/S_{n+1}^2 \leq \delta \leq 1/S_n^2$. 
Only the case $1/(2D_n^2Q_n) \leq \delta \leq 1/S_n^2$ is relevant, since the set in \eqref{time sequence} is empty otherwise. 
Then, 
\begin{equation*} 
\delta^r \#\{n : \delta < t_n < 2 \delta \} \le  \delta^r \, \cdot \delta\,  D^2 Q \, \cdot Q
\leq \frac{(DQ)^2}{S^{2(1+r)}}.
\end{equation*}
Hence, 
\begin{equation}\label{sequence_con}
DQ \lesssim S^{1+r} \qquad \Longrightarrow \qquad (t_n) \in \ell^{r,\infty}. 
\end{equation}

\subsection{Proof of Theorem ~\ref{main theorem} when $\boldsymbol{\alpha = n}$}
For the sequence $(t_n) \in \ell^{r,\infty}$ that we built, 
by \eqref{Pointwise_Lower_bound} and \eqref{Measure_X_R_Lebesgue} we get
\begin{equation}\label{maximal} 
\frac{\| \sup_{n}|e^{it_n\Delta} f_R|   \|_{L^2(0,1)}}{\|f_R\|_2} 
\gtrsim {S^{1/2}}\Big( \frac R{DQ}\Big)^{(n-1)/2}|X_R|^{1/2}
\simeq \Big( \frac{R^n}{S(DQ)^{n-1}} \Big)^{1/2}.                                                                                     
\end{equation}
Maximizing this quantity with the restrictions \eqref{measure_con} and \eqref{sequence_con}, 
we get
\begin{equation}\label{S_con}
S = R^{\frac{n}{r(n+1)+n} }, \quad DQ = S^{1+r} 
\quad \Longrightarrow \quad 
\frac{\| \sup_{n}|e^{it_n\Delta} f_R|   \|_{L^2(0,1)}}{\|f_R\|_2} \gtrsim R^{\frac{rn}{r(n+1)+n} }, 
\end{equation}
which disproves the maximal estimate when $s = rn/(r(n+1) + n) - \epsilon$ for arbitrary $\epsilon > 0$. 
Observe that the choice $S \geq R^{1/2}$ makes this argument be valid only for
\begin{equation}
R^{\frac{n}{r(n+1)+n}} \geq R^{1/2}
\qquad \Longleftrightarrow \qquad r \leq \frac{n}{n+1}. 
\end{equation}

\subsection{Heuristics for the fractal case $\boldsymbol{\alpha < n}$}
\label{sec:FractalHeuristics}
Let $\mu$ be an $\alpha$-dimensional measure. 
To disprove the fractal maximal estimate, 
we bound 
like in \eqref{maximal}, except now
\begin{align}\label{maximal_fractal}
\frac{\| \sup_{n}|e^{it_n\Delta} f_R|   \|_{L^2(d\mu)}}{\|f_R\|_2} &\gtrsim {S^{1/2}}\Big( \frac R{DQ}\Big)^{(n-1)/2} \mu(X_R)^{1/2}.                                                                                
\end{align}
To make $\mu(X_R)$ the largest possible, 
a heuristic argument
is to think of $\mu$ as the Hausdorff measure $\mathcal H^\alpha$. 
Then, recalling $X_R \subset [0,R/S^2] \times [0,1]^{n-1}$, 
we make $X_R$ scale like $[0,R/S^2] \times [0,1]^{n-1}$ when measured with $\mathcal H^\alpha$.
First, mimicking \textbf{Heuristic Argument 2} in Section~\ref{sec:Lebesgue_Measure}, 
 assume that the balls in $X_R$ are disjoint so that
\begin{align}
\label{FractalScalingMain}
\begin{split}
\mathcal H^{\alpha} (X_R) &\simeq Q\cdot \frac{D^2Q}{S^2}\cdot(DQ)^{n-1}\cdot \mathcal H^{\alpha} \Big( B_1\Big(\frac1S\Big)\times B_{n-1}\Big(\frac1R\Big) \Big) \\
                                            &\simeq \frac{(DQ)^{n+1}}{S^2} \frac RS \cdot \mathcal H^{\alpha} \Big(B_n\Big(\frac1R\Big)\Big) \\
                                            &\simeq \frac{(DQ)^{n+1}}{S^2} \frac RS \cdot \frac1{R^{\alpha}} \cdot \mathcal H^{\alpha} (B_n(1)), 
\end{split}
\end{align}
where we split the \textit{rectangle} $ B_1(1/S)\times B_{n-1}(1/R)$ into $R/S$ copies of $B_n(1/R)$. 
For the rectangle $[0,R/S^2] \times [0,1]^{n-1}$ we argue in two ways:
\begin{itemize}
	\item[\textbf{Scaling 1}] As in \eqref{FractalScalingMain}, split the rectangle in small cubes of size $R/S^2$:
	\begin{align}
	\label{FractalScaling1}   
	\begin{split}
	\mathcal H^{\alpha} ([0,R/S^2] \times [0,1]^{n-1}) 
	& \simeq \Big(\frac{S^2}R\Big)^{n-1}\cdot \mathcal H^{\alpha} \Big(B_n(R/S^2)\Big) \\
    &\simeq  \Big(\frac{S^2}R\Big)^{n-\alpha-1}\cdot \mathcal H^{\alpha} (B_n(1)).   
        \end{split}                                                        
	\end{align}
	\item[\textbf{Scaling 2}] Complete the rectangle to $[0,1]^n$ to get the same scaling as for the Lebesgue measure:
	\begin{equation}\label{FractalScaling2}
	\mathcal H^{\alpha} ([0,R/S^2] \times [0,1]^{n-1})  \simeq \frac{R}{S^2} \,   \mathcal H^{\alpha} ( [0,1]^n). 
	\end{equation}
\end{itemize}
\textbf{Scaling 1} will give a better result, but it will only work when $\alpha \geq n-1$. 
When $n/2 \leq \alpha \leq n - 1$ we will use \textbf{Scaling 2}.  

\begin{rk}
That \textbf{Scaling 1} works only for $\alpha \geq n-1$ is related to the fact that $(S^2/R)^{n - \alpha - 1} \to \infty$ when $\alpha < n-1$. 
Ultimately, this is a consequence of 
\begin{equation}
X_R \subset [0,R/S^2] \times [0,1]^{n-1} \to \{0\} \times [0,1]^{n-1} \quad \text{ as } R \to \infty,
\end{equation}
as mentioned in Remark~\ref{rk:Small_Measure}. 
\end{rk}

We compute both cases by equaling \eqref{FractalScalingMain} to either \eqref{FractalScaling1} or \eqref{FractalScaling2}. 

%
%
%

\subsubsection{Computations for Scaling 1}\label{sec:Scaling_1}
We equate \eqref{FractalScalingMain} and \eqref{FractalScaling1} to get 
\begin{equation}\label{fractal_con}
\frac{(DQ)^{n+1}}{SR^{\alpha}} \frac{R}{S^2} =\Big(\frac{S^2}R\Big)^{n-\alpha - 1}
\qquad \Longleftrightarrow \qquad 
\frac{(DQ)^{n+1}}{SR^n} =\Big(\frac{S}R\Big)^{2(n-\alpha)}.
\end{equation}
This restriction with the unchanged $\ell^{r,\infty}$ condition in \eqref{sequence_con} determines 
\begin{equation}\label{fractal_con2}
S = R^{\frac{2\alpha - n}{r(n+1) + 2\alpha - n}}, \qquad DQ = S^{1+r}.
\end{equation}
Hence, plugging the heuristic \eqref{FractalScalingMain} into \eqref{maximal_fractal}
and applying \eqref{fractal_con2}, 
we get 
\begin{align}\label{Heuristic_Result}
\frac{\| \sup_{n}|e^{it_n\Delta} f_R|   \|_{L^2(d\mu)}}{\|f_R\|_2} &\gtrsim   R^{\frac{n-\alpha}2} \, S^r
\simeq R^{\frac{n-\alpha}2 + \frac{r(2\alpha - n)}{r(n+1) + 2\alpha - n}},
\end{align}
which is the exponent in Theorem~\ref{main theorem} \eqref{MainTheorem_Item1}. 
Moreover, $S \geq R^{1/2}$ 
and \eqref{fractal_con2} require 
\begin{equation}\label{Condition_on_r}
\frac{2\alpha - n}{r(n+1) + 2\alpha - n} \geq \frac12 \quad \Longleftrightarrow \quad   r \leq \frac{2\alpha - n}{n+1}. 
\end{equation}

\subsubsection{Computations for Scaling 2}\label{sec:Scaling_2}
In this case we equate \eqref{FractalScalingMain} and \eqref{FractalScaling2} to get 
\begin{equation}\label{fractal_con_Scaling2}
\frac{(DQ)^{n+1}}{SR^{\alpha}} \frac{R}{S^2} = \frac{R}{S^2}
\qquad \Longleftrightarrow \qquad 
(DQ)^{n+1} = SR^\alpha. 
\end{equation}
Again, the unchanged $\ell^{r,\infty}$ condition in \eqref{sequence_con} with \eqref{fractal_con_Scaling2} determines 
\begin{equation}\label{fractal_con2_Scaling_2}
S = R^{\frac{\alpha}{r(n+1) + n}}, \qquad DQ = S^{1+r}.
\end{equation}
Plugging the heuristic \eqref{FractalScalingMain} into \eqref{maximal_fractal}, 
and applying \eqref{fractal_con_Scaling2},
we get 
\begin{align}\label{Heuristic_Result_Scaling2}
\frac{\| \sup_{n}|e^{it_n\Delta} f_R|   \|_{L^2(d\mu)}}{\|f_R\|_2} &\gtrsim   R^{\frac{n-\alpha}2} \, S^r 
\simeq R^{\frac{n-\alpha}2 + \frac{r \, \alpha}{r(n+1) +  n}},
\end{align}
which is the exponent in Theorem~\ref{main theorem} \eqref{MainTheorem_Item2}. 
The condition $S \geq R^{1/2}$ with \eqref{fractal_con2_Scaling_2} requires the same condition
\begin{equation}
\frac{\alpha}{r(n+1) + n} \geq \frac12 \quad \Longleftrightarrow \quad   r \leq \frac{2\alpha - n}{n+1}. 
\end{equation}

\subsection{Preparing the rigorous proof of Theorem~\ref{main theorem}}
To make the heuristics above rigorous, we will define a measure $\mu_R$ supported on $X_R$. 
Following the ideas in \cite{EP_V2}, 
we first redefine the parameters $S,D,Q$ by
\begin{equation}\label{Parameters}
R^A = DQ, \qquad R^B = \frac{D^2 Q}{R}, \qquad R^\sigma = S,
\end{equation}
which is unmade by
\begin{equation}
D = R^{1+B-A}, \qquad Q = R^{2A-B-1}.
\end{equation}
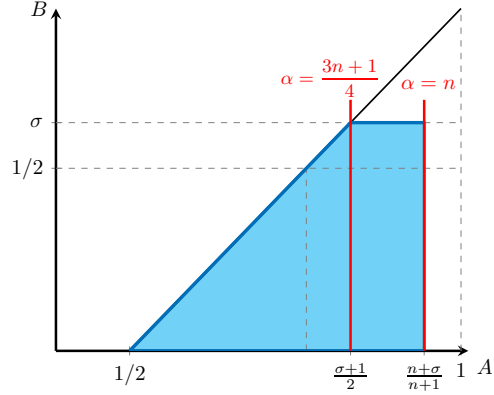
\begin{figure}[h]
\begin{center}
\begin{tikzpicture}[scale=0.8]
\begin{axis}[axis x line=center, axis y line=center, very thick,
xmin=0, xmax=11.2, ymin=0, ymax=4.5, xlabel={$A$}, ylabel={$B$}, xtick={0,2,8,10,11}, xticklabels={$0$, $1/2$, $\frac{\sigma+1}{2}$,$\frac{n+\sigma}{n+1}$, $1$}, ytick={2.4,3}, yticklabels={$1/2$,$\sigma$},   
xlabel style={below right},
ylabel style={left}
]
\filldraw[fill=ProcessBlue!50] (2,0)--(10,0)--(10,3)--(8,3)--(2,0);
\addplot[black, domain=0:11, thick] {x/2-1};
\addplot[gray, dashed, thin] coordinates {(0,3) (11,3)};
\addplot[gray, dashed, thin] coordinates {(0,2.4) (11,2.4)};
\addplot[gray, dashed, thin] coordinates {(6 .8,0) (6.8,2.4)};
\addplot[gray, dashed, thin] coordinates {(11,0) (11,4.5)};
\addplot[NavyBlue, domain=2:8, ultra thick] {x/2-1};
\addplot[NavyBlue, domain=8:10, ultra thick] {3};
\addplot[NavyBlue, domain=2:10, ultra thick] {0};
\addplot[red, very thick] coordinates {(8,0) (8,3.3)};
\addplot[red, very thick] coordinates {(10,0) (10,3.3)};
\node[above=5.8cm, right=3.6cm, red] {\small{$\alpha= \displaystyle{\frac{3n+1}{4}}$}};
\node[above=5.7cm, right=5.6cm, red] {$\alpha=n$};
\end{axis}
\end{tikzpicture}
\caption{Region for the parameters $A,B$ defined in \eqref{Parameters} for a fixed $\sigma \in  [1/2,1]$.}
\end{center}
\end{figure}
The conditions $S \geq R^{1/2}$ and the elementary $Q \geq 1$, together with forcing the rectangles in each $q \simeq Q$ level of $X_R$ to be disjoint, 
impose the basic restrictions 
\begin{equation}\label{Basic_Conditions}
\frac12 < \sigma \le 1, \qquad 0 \leq A \leq 1, \qquad 0 \leq B \leq \sigma, 
\qquad 2A \geq B + 1.
\end{equation}

\subsection{Proof using Scaling 1} 
From \eqref{fractal_con} and \eqref{fractal_con2} we set
\begin{equation}\label{fractal_con_Scaling1_Bis}
\frac{(DQ)^{n+1}}{SR^n} =\Big(\frac{S}R\Big)^{2(n-\alpha)}, 
\qquad  
DQ = S^{1+r}, 
\end{equation}
This is rewritten in terms of  $A$ and $\sigma$ as
\begin{equation}\label{fractal_con_Exponents}
A(n+1) = 2\alpha - n + 2\sigma(n-\alpha + 1/2), 
\qquad A = \sigma (1+r), 
\end{equation}
which are independent\footnote{Moreover, the heuristic Sobolev exponent in \eqref{Heuristic_Result} depends on $\sigma$ (hence on $A$), but not on $B$.}
of $B$. 
The first condition in \eqref{fractal_con_Exponents} can be rewritten as 
\begin{align}
A(n+1) & = \alpha + \sigma + (2\sigma-1)(n-\alpha)
\label{fractal_con_Exponents_2} \\
& = n + \sigma - 2(1-\sigma)(n-\alpha).
\label{fractal_con_Exponents_3}
\end{align}
Based on \eqref{FractalScaling1}, define the measure 
\begin{equation}\label{measure}
\mu_R(E) = \frac{|E\cap X_R|}{|X_R|} \, \Big(\frac{S^2}{R}\Big)^{n - 1 - \alpha}
\end{equation}
so that $\mu(X_R) = (S^2/R)^{n-1-\alpha}$.
Since the computations in Section~\ref{sec:Scaling_1} give the desired Sobolev exponents, 
it suffices to prove that $\mu_R$ is an $\alpha$ dimensional measure. 
We first simplify the expression for $\mu_R$. 
From \eqref{X_R_in_terms_of_tildeXR}, $|X_R| = |\widetilde X_R| R/S^2$. 
Now, 
\begin{equation}\label{Measure_Of_X_R_Calculation}
Q^{n+1} h_1(Q) h_2(Q)^{n-1} = \frac{(DQ)^{n+1}}{ S R^n }
=   \left( \frac{S}{R} \right)^{2(n-\alpha)} \ll 1
\end{equation}
because of \eqref{fractal_con_Scaling1_Bis}, 
so Lemma~\ref{Lemma_AnChuPierce} implies
\begin{equation}\label{Measure_Of_X_R}
 |\widetilde X_R| 
 \simeq Q^{n+1} h_1(Q) h_2(Q)^{n-1} 
 \simeq  \left( \frac{S}{R} \right)^{2(n-\alpha)}
 \quad \Longrightarrow \quad
 |X_R| \simeq \frac{R}{S^2} \, \left( \frac{S}{R} \right)^{2(n-\alpha)}. 
\end{equation}
Therefore, we rewrite the measure $\mu_R$ as
\begin{equation}\label{Measure_In_Compact_Form}
\frac{1}{|X_R|} \, \Big(\frac{S^2}{R}\Big)^{n - 1 - \alpha}
\simeq   R^{n - \alpha} 
\quad \Longrightarrow \quad 
\mu_R(E) \simeq |E\cap X_R| \, R^{n-\alpha}.
\end{equation}

We now prove that $\mu_R$ is $\alpha$-dimensional.
Let $r > 0$ and $B_r$ any ball of radius $r$. 
We want to prove that $\mu_R(B_r) \lesssim r^\alpha$ with a constant independent of $R$. 
We will exploit the freedom for $B$ by choosing it as large as possible. 
Direct computation from \eqref{fractal_con_Exponents_3}
shows that there are two different ranges, which we treat separately:
\begin{enumerate}
	\item $A \leq (1+\sigma)/2 \Longleftrightarrow \alpha \leq (3n+1)/4$. In this case, we choose $B = 2A - 1$, which corresponds to $Q = 1$.
	\item $A \geq (1+\sigma)/2 \Longleftrightarrow \alpha \geq (3n+1)/4$. In this case, we choose $B = \sigma$. 
\end{enumerate}

\subsubsection{When $B = 2A - 1$ and $Q=1$.}\label{sec:Scaling1_Q_1}
In this case, 
the set $X_R$ becomes particularly easy because there is no union in $q \simeq Q = 1$, so 
\begin{equation}\label{baseset_Qequals}
X_R = \bigcup_p B_1 \Big( \frac R{D^2} p_1, \frac cS \Big) \times B_{n-1}\Big(\frac1D p', \frac cR \Big) \subset [0,R/S^2] \times [0,1]^{n-1}.
\end{equation}
In particular, all rectangles are disjoint. 
We measure $B_r$ depending on $r$:  
\begin{itemize}
	\item $r \leq 1/R$. 
	Bound $|B_r \cap X_R| \leq |B_r| = r^n$ so
	\begin{equation}
	\mu_R(B_r) \leq r^n\, R^{n-\alpha} = r^\alpha (rR)^{n-\alpha} \leq r^\alpha. 
	\end{equation}
	
	\item $1/R < r \leq 1/D$. The ball $B_r$ intersects at most one rectangle, so 
	\begin{equation}
	\mu_R(B_r) \leq \frac{r}{R^{n-1} } \, R^{n-\alpha}  
	= r^\alpha \Big( \frac{1}{rR} \Big)^{\alpha - 1}	\leq r^{\alpha}, \qquad \text{ if } \alpha \geq 1.
	\end{equation}	 
	
\item $1/D < r \leq 1/S$. (Observe that $S \leq D = DQ = S^{1+r}$). 
Now $B_r$ intersects at most $(rD)^{n-1}$ rectangles in direction $x'$ and one in direction $x_1$, each of which can contribute at most with a measure $r/R^{n-1}$,
so
\begin{equation}
\mu_R(B_r) \leq (rD)^{n-1} \, \frac{r}{R^{n-1}} \,  R^{n-\alpha}
= r^n\, \frac{D^{n+1}}{R^{\alpha}} \frac{R}{D^2}
= r^n R^{A(n+1) - \alpha -  B }.
\end{equation}
By bounding $r \le R^{-\sigma}$ and using the condition \eqref{fractal_con_Exponents_2}, we get
\begin{equation}
\mu_R(B_r) \leq 
r^\alpha\,  R^{\sigma  - B - (1 - \sigma)(n-\alpha)}
\leq r^\alpha,
\end{equation}
because it follows from \eqref{fractal_con_Exponents_3}, $2A=B+1$ and $B \leq \sigma \leq 1$ that
\begin{equation}\label{Intermediate_Balls_Condition}
\begin{split}
\sigma - B -(1-\sigma)(n-\alpha)
& =\sigma - B + \frac{B+1}{2}\frac{n+1}{2} - \frac{\sigma}{2} - \frac{n}{2} \\
& =  \frac{\sigma}{2} + B \frac{n-3}{4} - \frac{n-1}{4} \\
& \leq (\sigma - 1) \frac{n-1}{4} \leq 0.
\end{split}
\end{equation}
	
	\item $1/S < r \leq R/D^2 = 1/R^B$. 
	The ball intersects at most $(rD)^{n-1}$ rectangles in direction $x'$ and one in direction $x_1$, 
	each contributing with $1/(SR^{n-1})$. 
	Then, by \eqref{fractal_con_Scaling1_Bis}, 
	\begin{equation}
	\mu_R(B_r) \leq \frac{(rD)^{n-1}}{SR^{n-1} }R^{n-\alpha}
	= r^{n-1} \frac{D^{n+1}}{SR^\alpha } \frac{R}{D^2} 
	= r^{n-1} \, \Big( \frac{S^2}{R} \Big)^{n-\alpha}  \,  \frac{1}{R^B}. 
\end{equation}	 
If $\alpha \geq n-1$, we bound $r^{n-1} = r^{\alpha} \, r^{n - 1 - \alpha} \leq r^\alpha \, S^{1 - (n-\alpha)}$,
so
\begin{equation}
\mu_R(B_r) \leq r^\alpha \,\frac{S}{R^B} \,   \Big( \frac{S}{R} \Big)^{n-\alpha}
=  r^\alpha \, R^{ \sigma - B - (1- \sigma)(n-\alpha) }
\leq r^\alpha,
\end{equation}
where the last inequality follows from \eqref{Intermediate_Balls_Condition}.

\item $1/R^B < r \leq R/S^2$. The ball intersects at most $(rD)^{n-1}$ rectangles in direction $x'$ and $R^B r$ in direction $x_1$, each with a contribution of $1/(SR^{n-1})$. 
	Then, by \eqref{fractal_con_Scaling1_Bis}, 
	\begin{equation}
	\mu_R(B_r) \leq r^n \frac{D^{n-1}R^B}{SR^{n-1}} \, R^{n-\alpha}
	= r^n \frac{D^{n+1}}{SR^n} \, R^{n-\alpha}
	= r^\alpha \Big(  r \frac{S^2}{R} \Big)^{n-\alpha}
	\leq r^\alpha.
\end{equation}

\item $R/S^2 < r \leq 1$. The ball intersects at most $(rD)^{n-1}$ rectangles in direction $x'$ and $D^2 /S^2$ in direction $x_1$.
So if $\alpha \geq n-1$, from \eqref{fractal_con_Scaling1_Bis} we get
\begin{equation}
\mu_R(B_r) \leq r^{n-1} \frac{D^{n+1}}{S^2} \, \frac{1}{SR^{n-1}} \,  R^{n-\alpha}
= r^{n - 1} \, \frac{R}{S^2} \Big( \frac{S^2}{R} \Big)^{n-\alpha}
= r^\alpha \Big( r \frac{S^2}{R} \Big)^{n - 1- \alpha}
\leq r^\alpha.
\end{equation}
	
\end{itemize}
Hence $\mu_R$ is $\alpha$-dimensional, with underlying constant $\langle \mu_R \rangle_\alpha \simeq 1$.

\subsubsection{When $B= \sigma$ and  $Q\gg1$}\label{sec:Scaling1_Q_Large}
Now $X_R$ is a superposition of grids of rectangles, one for each $q \simeq Q$. 
However, by the choice of $B=\sigma$, for each fixed $x'$ the rectangles form a single tube of length $R/S^2$ in direction $x_1$. 
To control intersections in direction $x'$, we slightly adapt a lemma proved by Lucà and Ponce-Vanegas in \cite{LP-V}.
\begin{lem}[Lemma 13 in \cite{LP-V}]
\label{Lemma_LPV}
Let $X_R = \cup_{q} \cup_p B_1 \times B_{n-1} \subset [0,R/S^2] \times [0,1]^{n-1}$ as defined in \eqref{baseset}, 
which is a union of $(DQ)^{n+1}/S^2$ rectangles.
Let $\epsilon > 0$. Then, for $R \gg_\epsilon 1$, at least $ (DQ)^{n+1}/(S^2Q^\epsilon)$ such rectangles have the following property: 
if $x=(x_1,x')$ and $y=(y_1,y')$ are centers of two such rectangles, 
\begin{equation}
 |x'-y'| > \frac1{DQ^{\frac n{n-1}}}.
\end{equation}
\end{lem}
Call $Y_R = \bigcup_q \bigcup_p B_1 \times B_{n-1} \subset X_R$ 
the union restricted to the subcollection of rectangles given by Lemma~\ref{Lemma_LPV}.
Instead of working with $\mu_R$ as in \eqref{measure}, 
we define 
\begin{equation}\label{measure_Y}
\nu_R(E) = \frac{|E\cap Y_R|}{|Y_R|} \, R^{-\epsilon} \,  \Big( \frac{S^2}{R}  \Big)^{n-1-\alpha}
\end{equation}
for some $\epsilon > 0$. 
Observe that the rectangles in $Y_R$ are pairwise disjoint because 
\begin{equation}\label{Disjointness_Property}
 \frac{1}{R} \leq \frac{1}{D Q^{n/(n-1)}} \quad \Longleftrightarrow \quad (n+1)A \leq n + B, 
\end{equation} 
which holds in view of $B=\sigma$ and \eqref{fractal_con_Exponents_3}. 
Hence, comparing with $|X_R|$ in \eqref{Measure_Of_X_R},
\begin{equation}\label{Measure_Y_R}
|Y_R| \simeq \frac{1}{R^\epsilon} \, \frac{(DQ)^{n+1}}{S^2} \frac{1}{SR^{n-1}} \simeq \frac{1}{R^\epsilon}|X_R|.
\end{equation}
Together with \eqref{Measure_In_Compact_Form}, we get $\nu_R(E) \simeq  |E \cap Y_R| \, R^{n-\alpha}$. 
We now prove that $\nu_R$ is an $\alpha$-dimensional measure.
Let $B_r$ be a ball of radius $r$. 
Again, recall that $B=\sigma$ implies that for each fixed $x'$, the rectangles form a single tube in direction $x_1$. 
\begin{itemize}
	\item If $r \leq 1/R$, then $|Y_R \cap B_r| \leq |B_r| = r^n$, so
	\begin{equation}
	\nu_R(B_r) \leq r^n R^{n-\alpha} = r^\alpha (rR)^{n-\alpha} \leq r^\alpha. 
	\end{equation}
	
	\item $1/R < r \leq 1/(DQ^{n/(n-1)})$. Then $B_r$ intersects at most one tube, for which 
	\begin{equation}
	\nu_R(B_r) \leq  \frac{r}{R^{n-1}} \,  R^{n-\alpha} = r^\alpha \Big(\frac{1}{rR} \Big)^{\alpha - 1} \leq r^\alpha, \qquad \text{ if } \alpha \geq 1. 
	\end{equation}
	
	\item If $1/(DQ^{n/(n-1)}) < r \leq R/S^2$, then $B_r$ intersects at most $(rDQ^{n/(n-1)})^{n-1}$ tubes, 
	so 
	\begin{align}
	\nu_R(B_r)  \leq  r^{n-1}  D^{n-1}  Q^n \, \frac{r}{R^{n-1}}\,   R^{n-\alpha}
	= r^n \, \frac{D^{n-1}Q^n}{R^{\alpha - 1}} 
	= r^n \, \frac{R^{A(n+1)}}{R^\alpha \, R^B}.
	\end{align}
	Since $B=\sigma$, by \eqref{fractal_con_Exponents_2} we get $A(n+1) - \alpha - B = (2\sigma - 1) (n-\alpha)$.
	Together with $r \leq R/S^2 = R^{1 - 2\sigma}$, this implies
	\begin{equation}
	\nu_R(B_r) \leq  r^n \, R^{(2\sigma - 1)(n-\alpha)} = r^\alpha (rR^{2\sigma - 1})^{n-\alpha} \leq r^\alpha. 
	\end{equation}
	
	\item If $R/S^2 < r  \leq 1$, then $B_r$ intersects again at most $(rDQ^{n/(n-1)})^{n-1}$ tubes, but the measure of each tube is now $(R/S^2) (1/R^{n-1})$, so 
	 \begin{align}
	\nu_R(B_r) 
	\leq r^{n-1}  D^{n-1}  Q^n \, \frac{R}{S^2}\frac{1}{R^{n-1}}\,   R^{n-\alpha} 
    = r^{n-1} \frac{R}{S^2} \frac{R^{A(n+1)}}{R^\alpha\, R^B} 
	\end{align}
	Again, by $B=\sigma$ and \eqref{fractal_con_Exponents_2} we get $A(n+1) - \alpha - B = (2\sigma - 1)(n- \alpha)$, so
	\begin{equation}
	\nu_R(B_r) \leq r^{n-1} \,  R^{1-2\sigma} \, R^{(2\sigma-1)(n- \alpha)} = r^\alpha \, (rR^{2\sigma-1})^{n - \alpha - 1} \leq r^\alpha, 
	\end{equation}
	if $\alpha \geq n-1$. 
\end{itemize}
Hence, $\nu_R$ is $\alpha$-dimensional, with constant $ \langle \nu_R \rangle_\alpha \simeq 1$. 

\subsubsection{Concluding the proof}
Let $n-1 \leq \alpha  < n$. Let $R \gg 1$. 
Define the measure $\lambda_R = \mu_R$ if $\alpha \leq (3n+1)/4$, 
and $\lambda_R = \nu_R$  if $\alpha \geq (3n+1)/4$.
In both cases, $\lambda_R$ is an $\alpha$-dimensional measure and $\langle \lambda_R \rangle_\alpha \simeq 1$. 
By \eqref{Pointwise_Lower_bound}, \eqref{measure} and \eqref{measure_Y}, 
\begin{equation}
\frac{\| \sup_{n}|e^{it_n\Delta} f_R|   \|_{L^2(d\lambda_R)}}{\|f_R\|_2} 
\gtrsim S^{\frac12} \, \Big( \frac R{DQ}\Big)^{\frac{n-1}2} \lambda_R(X_R)^{\frac12}
\simeq   S^{\frac12} \, \Big( \frac R{DQ}\Big)^{\frac{n-1}2} \, R^{-\epsilon} \,  \Big(\frac{S^2}{R}  \Big)^{\frac{n-1 - \alpha}{2}}.                                  
\end{equation}
Solving for $DQ$ and $S$ from \eqref{fractal_con_Scaling1_Bis}, we get $S = R^{\frac{2\alpha - n}{r(n+1) + 2\alpha - n}}$. 
Since $S \geq R^{1/2}$, this forces $r \leq (2\alpha - n)/(n+1)$. 
Hence 
\begin{equation}
\frac{\| \sup_{n}|e^{it_n\Delta} f_R|   \|_{L^2(d\lambda_R)}}{\|f_R\|_2} 
\gtrsim R^{\frac{n-\alpha}{2} + \frac{r(2\alpha - n)}{r(n+1) + 2\alpha - n} - \epsilon}.                     
\end{equation}
If we assume the fractal maximal estimate for $\alpha$-dimensional measures in $H^s$, then 
\begin{equation}
R^{\frac{n-\alpha}{2} + \frac{r(2\alpha - n)}{r(n+1) + 2\alpha - n} - \epsilon }
\lesssim \langle \lambda_R \rangle_\alpha \, R^s \lesssim R^s \qquad \forall R \gg_\epsilon 1,
\end{equation}
which gives a contradiction for $s < \frac{n-\alpha}{2} + \frac{r(2\alpha - n)}{r(n+1) + 2\alpha - n}$. 
Finally, when $r = \frac{2\alpha - n}{n+1}$, 
the exponent is $\frac{n-\alpha}{2}  + \frac{2\alpha - n}{2(n+1)}$,  
which immediately implies 
$s_c^{\operatorname{max}}(\alpha, \ell^{r,\infty}) \geq \frac{n-\alpha}{2}  + \frac{2\alpha - n}{2(n+1)}$
for $r \geq \frac{2\alpha - n}{n+1}$.
Combining these two gives the statement in Theorem~\ref{main theorem}.\eqref{MainTheorem_Item1}. 
$\qed$

\subsection{Proof using Scaling 2}
From \eqref{fractal_con_Scaling2} and \eqref{fractal_con2_Scaling_2}, 
now set 
\begin{equation}\label{Conditions_Scaling2}
(DQ)^{n+1} = SR^\alpha, \qquad DQ = S^{1+r}, 
\end{equation}
which is rewritten in terms of $A$ and $\sigma$ as 
\begin{equation}\label{Conditions_Exponents_Scaling2}
A(n+1) = \sigma + \alpha, \qquad A = \sigma(1+r). 
\end{equation}
Define the measure
\begin{equation}
\mu_R(E) = \frac{|E \cap X_R|}{|X_R|} \, \frac{R}{S^2}. 
\end{equation}
Proceeding like in \eqref{Measure_Of_X_R_Calculation} and \eqref{Measure_Of_X_R}, 
we get 
\begin{equation}\label{Measure_X_R_Scaling2}
|X_R| \simeq \frac{(DQ)^{n+1}}{SR^n} \,  \frac{R}{S^2} 
 \simeq \frac{1}{R^{n-\alpha}} \, \frac{R}{S^2}
\quad \Longrightarrow \quad 
\mu_R(E) \simeq |E \cap X_R| \, R^{n-\alpha}. 
\end{equation}
We prove now that $\mu_R$ is an $\alpha$-dimensional measure. 
We separate two cases again:

\subsubsection{When $B = 2A - 1$ and $Q = 1$.}
We repeat the procedure in Section~\ref{sec:Scaling1_Q_1}
by working with the simpler, disjoint structure of $X_R$ in \eqref{baseset_Qequals}. 
Let $r > 0$. 
\begin{itemize}

	\item If $r \leq 1/R$, 
	bound $|B_r \cap X_R| \leq |B_r| = r^n$ so
	\begin{equation}
	\mu_R(B_r) \leq r^n\, R^{n-\alpha} = r^\alpha (rR)^{n-\alpha} \leq r^\alpha. 
	\end{equation}
	
	\item If $1/R < r \leq 1/D$, then $B_r$ intersects at most one rectangle, so 
	\begin{equation}
	\mu_R(B_r) \leq  \frac{r}{R^{n-1}} \, R^{n-\alpha}
	= r^\alpha \Big( \frac{1}{rR} \Big)^{\alpha - 1}	\leq r^{\alpha}, \quad \text{ if } \alpha - 1 \geq 0.
	\end{equation}	 
	
	\item If $1/D < r \leq 1/S$, $B_r$ intersects at most $(rD)^{n-1}$ rectangles in direction $x'$ and one in direction $x_1$, 
	each contributing at most with $r/R^{n-1}$. 
	By \eqref{Conditions_Scaling2}, 
	\begin{equation}
	\mu_R(B_r) \leq (rD)^{n-1} \,  \frac{r}{R^{n-1}} \, R^{n-\alpha}
	= r^n\, \frac{D^{n+1}}{R^{\alpha}} \frac{R}{D^2}
	=  r^n\, \frac{S}{R^B}.
	\end{equation}
	Then, since $r^{n-\alpha} \le 1/S^{n-\alpha}$, 
	\begin{equation}
	\mu_R(B_r) \leq 
	\frac{r^\alpha}{R^B \, S^{n-\alpha - 1}} 
	\leq r^\alpha, 
	\qquad \text{ if } \alpha \leq n-1.
	\end{equation}

	\item If $1/S < r \leq 1/R^B$, $B_r$ intersects at most $(rD)^{n-1}$ rectangles in direction $x'$ and one in direction $x_1$, 
	each contribution being $1/(SR^{n-1})$. 
	By \eqref{Conditions_Scaling2}, 
	\begin{equation}
	\mu_R(B_r) \leq \frac{(rD)^{n-1}}{SR^{n-1}}R^{n-\alpha}
	= r^{n-1} \,  \frac{D^{n+1}}{SR^\alpha}\,  \frac{R}{D^2}
	= r^\alpha \, \frac{r^{n-1-\alpha}}{R^B}.
	\end{equation}	 
	Since $r \leq 1/R^B$, we get 
	\begin{equation}
	\mu_R(B_r) \leq   \frac{r^\alpha}{R^B \ R^{B(n-\alpha - 1)}} = \frac{r^\alpha}{R^{B(n-\alpha)}} \leq r^\alpha, 
	\qquad \text{ if } \alpha \leq n-1.
	\end{equation}
	
	\item  $1/R^B < r \leq 1$. The ball intersects at most $(rD)^{n-1}$ rectangles in direction $x'$ and $rR^B$ in direction $x_1$, 
	each contribution being $1/(SR^{n-1})$. 
	By \eqref{Conditions_Scaling2}, 
	\begin{equation}
	\mu_R(B_r) \leq \frac{(rD)^{n-1} rR^B}{SR^{n-1}}R^{n-\alpha}
	= r^n \,  \frac{D^{n+1}}{SR^\alpha}\,  \frac{R}{D^2} \, R^B
	= r^n 
	\leq r^\alpha.
	\end{equation}	 

\end{itemize}
Hence, $\mu_R$ is an $\alpha$-dimensional measure with constant $\langle \mu_R \rangle \simeq 1$.

\subsubsection{When $B = \sigma$ and $Q \gg 1$.}

As in Section~\ref{sec:Scaling1_Q_Large}, 
using Lemma~\ref{Lemma_LPV} 
we extract from $X_R$ a subcollection $Y_R$ of $(DQ)^{n+1}/(S^2Q^\epsilon)$ rectangles. 
The condition for disjointness in the $x'$ direction in \eqref{Disjointness_Property} is also satisfied because of \eqref{Conditions_Scaling2}. 
As before, since $B=\sigma$, for each fixed $x'$ the rectangles form a single tube of length $R/S^2$ in direction $x_1$. 
Define
\begin{equation}
\nu_R(E) = \frac{|E \cap Y_R|}{|Y_R|} \, R^{-\epsilon} \,  \frac{R}{S^2} . 
\end{equation}
By disjointness, proceeding like in \eqref{Measure_Y_R}, we get $|Y_R| \simeq |X_R|/R^{\epsilon}$, 
so from \eqref{Measure_X_R_Scaling2} we get $\nu_R(E) \simeq |E \cap Y_R| R^{n-\alpha}$. 
Let now $r > 0$. 
\begin{itemize}

	\item If $r \leq 1/R$, we bound $|B_r \cap Y_R| \leq |B_r| = r^n$, so
	\begin{equation}
	\nu_R(B_r) \leq r^n R^{n-\alpha} = r^\alpha (rR)^{n-\alpha} \leq r^\alpha. 
	\end{equation}
	
	\item If $1/R < r \leq 1/(DQ^{n/(n-1)})$, $B_r$ intersects at most one tube, for which 
	\begin{equation}
	\nu_R(B_r) \leq  \frac{r}{R^{n-1}} \,  R^{n-\alpha} = r^\alpha \Big(\frac{1}{rR} \Big)^{\alpha - 1} \leq r^\alpha,
	\quad \text{ if } \alpha \geq 1.
	\end{equation}
	
	\item If $1/(DQ^{n/(n-1)}) < r \leq 1$, $B_r$ intersects at most $(rDQ^{n/(n-1)})^{n-1}$ tubes: 
	\begin{align}
	\nu_R(B_r)  \leq r^{n-1}D^{n-1}Q^n \, \frac{r}{R^{n-1}}\,   R^{n-\alpha}
	= r^n \, \frac{R^{A(n+1)}}{R^\alpha \, R^B}.
	\end{align}
	Since $B=\sigma$, from \eqref{Conditions_Exponents_Scaling2} we get $A(n+1) = \alpha + B$, so 
	\begin{equation}
	\nu_R(B_r) \leq  r^n \leq r^\alpha. 
	\end{equation}
	
\end{itemize}
Hence, $\nu_R$ is an $\alpha$-dimensional measure with constant $\langle \nu_R \rangle \simeq 1$.

\subsubsection{Concluding the proof}

Let $n/2\leq \alpha  < n$. Let $R \gg 1$. 
Define the measure $\lambda_R = \mu_R$ or $\lambda_R = \nu_R$, depending on the value of $\alpha$.
In both cases, $\lambda_R$ is an $\alpha$-dimensional measure and $\langle \lambda_R \rangle_\alpha \simeq 1$. 
By \eqref{Pointwise_Lower_bound}, 
\begin{equation}
\frac{\| \sup_{n}|e^{it_n\Delta} f_R|   \|_{L^2(d\lambda_R)}}{\|f_R\|_2} 
\gtrsim S^{\frac12} \, \Big( \frac R{DQ}\Big)^{\frac{n-1}2} \lambda_R(X_R)^{\frac12}
\simeq   S^{\frac12} \, \Big( \frac R{DQ}\Big)^{\frac{n-1}2} \, R^{-\epsilon} \, \Big( \frac{R}{S^2}\Big)^{1/2}.
\end{equation}
Solving for $DQ$ and $S$ from \eqref{Conditions_Scaling2} we get $S = R^{\frac{\alpha}{r(n+1) + n}}$ and $DQ = S^{1+r}$,
which forces $r \leq (2\alpha - n)/(n+1)$. 
Hence, 
\begin{equation}
\frac{\| \sup_{n}|e^{it_n\Delta} f_R|   \|_{L^2(d\lambda_R)}}{\|f_R\|_2} 
\gtrsim R^{\frac{n-\alpha}{2} + \frac{r\alpha}{r(n+1) + n} - \epsilon}.                     
\end{equation}
If we assume the fractal maximal estimate for $\alpha$-dimensional measures in $H^s$, then 
\begin{equation}
R^{\frac{n-\alpha}{2} + \frac{r\alpha}{r(n+1) + n} - \epsilon }
\lesssim \langle \lambda_R \rangle_\alpha \, R^s \lesssim R^s \qquad \forall R \gg_\epsilon 1,
\end{equation}
which gives a contradiction for $s < \frac{n-\alpha}{2} + \frac{r\alpha}{r(n+1) + n}$. 
Evaluating this exponent at $r = \frac{2\alpha - n}{n+1}$, 
we immediately get 
$s_c^{\operatorname{max}}(\alpha, \ell^{r,\infty}) \geq \frac{n-\alpha}{2}  + \frac{2\alpha - n}{2(n+1)}$
for $r \geq \frac{2\alpha - n}{n+1}$.
Combining these two gives the statement in Theorem~\ref{main theorem}.\eqref{MainTheorem_Item2}. 
$\qed$

\section{Proof of Theorem~\ref{main_theorem_1D} - Dimension $n=1$}

\subsection{When $\boldsymbol{\alpha = 1}$}\label{sec:1D_Lebesgue}
We reprove the result by Dimou and Seeger \cite{DS}. 
We work with the function $g_R$ defined in \eqref{Def_Of_ID}. 
In \eqref{con1}  we got 
\begin{equation}\label{con1_BIS}
 S|x_1+2Rt| \le c, \quad S^2t \le c
 \quad  \Longrightarrow \quad 
 |e^{it\Delta}g_R(x_1)| \simeq 1.
\end{equation} 
As we did for $n \geq 2$ in Section~\ref{sec:Proof_Main_Thm}, we work with $S \geq R^{1/2}$. 
Define
\begin{equation}\label{1D_Sets}
T_R = \Big\{ \, \frac{p}{Q} \, \,  : \, \,  1 \leq p \leq \frac{Q}{S^2}  \, \Big\}, 
\quad \text{ and } \quad 
X_R = \bigcup_{p=1}^{Q/S^2} B\Big(\frac{R}{Q}\, p, \frac1S \Big).
\end{equation} 
We impose the restrictions $S^2 \leq Q$ and $Q \leq RS$, 
the first one for the sets in \eqref{1D_Sets} not to be empty, and the second for the intervals in $X_R$ to be disjoint. 
Actually, since $X_R \subset [0,R/S^2]$, we get the maximum measure when
\begin{equation}\label{1D_Measure_Condition}
Q = RS \quad \Longrightarrow \quad |X_R| \simeq \frac{Q}{S^2} \, \frac{1}{S} = \frac{R}{S^2}.
\end{equation}
On the other hand, observe that $T_R \subset [1/Q, 1/S^2]$. 
Asking $R_n$ to be a dyadic sequence that grows fast enough, all intervals $[1/Q_n, 1/S_n^2]$ are disjoint. 
Moreover, for any $\delta > 0$, 
place it in $ 1/S_{n+1}^2 \leq \delta \leq 1/S_n^2$ so that
\begin{equation}\label{1D_Sequence_Condition_Proto}
\delta^r \, \#\{ \, t \in T_{R_n} \,  : \,  \delta < t < 2 \delta \,  \} \le  \delta^{1+r} Q_n \le \frac{Q_n}{S_n^{2(1+r)}}.
\end{equation} 
Hence, defining our sequence of times as $\{ t_n \} = \bigcup_n T_{R_n}$, 
\begin{equation}\label{1D_Sequence_Condition}
Q \lesssim S^{2(1+r)} \qquad \Longrightarrow \qquad \{ t_n \} \in \ell^{r,\infty}. 
\end{equation}
Gathering all, we get 
\begin{equation} \label{1D_L2_Bound}
\frac{\| \sup_{n}|e^{it_n\Delta} g_R|   \|_{L^2}}{\|g_R\|_2} \gtrsim  S^{1/2} \, |X_R|^{1/2} 
= \Big(\frac{R}{S}\Big)^{1/2}.                                                                                        
\end{equation}
Combining \eqref{1D_Measure_Condition} and taking the smallest $S$ in \eqref{1D_Sequence_Condition},
we get $S = R^{1/(1+2r)}$, hence 
\begin{equation}\label{1D_Exponent}
\frac{\| \sup_{n}|e^{it_n\Delta} g_R|   \|_{L^2}}{\|g_R\|_2} \gtrsim R^{\frac{r}{2r+1}}.                                                                                        
\end{equation}
Hence the maximal estimate cannot work in $H^s$ for $s < r/(2r+1)$. 
Observe also that the choice $S \geq R^{1/2}$ restricts $r \leq 1/2$. 

Since when $r=1/2$ the exponent in \eqref{1D_Exponent} is 1/4, 
this immediately implies that $s_c^{\operatorname{max}}(\ell^{r,\infty}) \geq 1/4$ for all $r \geq 1/2$. 
$\qed$

\subsection{When $\boldsymbol{\alpha  <  1}$}\label{sec:1D_Fractal}

With the same function $g_R$, the constructions \eqref{1D_Sets} 
and the condition \eqref{1D_Sequence_Condition} intact, 
we alter condition \eqref{1D_Measure_Condition}. 
Heuristically, we want $X_R$ and $[0,R/S^2]$ to have the same $\mathcal H^\alpha$-scaling. 
Since  
\begin{equation}
\mathcal H^\alpha (X_R) \simeq \frac{Q}{S^2} \frac{1}{S^\alpha} \, \mathcal H^\alpha ( [ 0,1]), 
\quad \text{ and } \quad 
\mathcal H^\alpha \Big(  \Big[ 0,\frac{R}{S^2} \Big] \Big) \simeq \Big( \frac{R}{S^2} \Big)^\alpha \,\mathcal  H^\alpha ( [ 0,1]),  
\end{equation}
we choose
\begin{equation}\label{1D_Fractal_Condition}
\frac{Q}{S^2} \frac{1}{S^\alpha}  = \Big( \frac{R}{S^2} \Big)^\alpha 
\quad \Longleftrightarrow \quad Q = R^\alpha \, S^{2-\alpha}.
\end{equation}
In particular, $Q = RS (S/R)^{1-\alpha} \ll RS$. 
With that, heuristically we expect to get  
\begin{equation} 
\frac{\| \sup_{n}|e^{it_n\Delta} g_R|   \|_{L^2(\mathcal H^\alpha)}}{\|g_R\|_2} \gtrsim  S^{1/2} \, \mathcal H^\alpha (X_R)^{1/2} 
\simeq S^{1/2} \, \Big(\frac{R}{S^2}\Big)^{\alpha/2}
=   \Big( \frac{R^\alpha}{S^{2\alpha - 1}}\Big)^{1/2}.                       
\end{equation}
Taking the smallest $S$ in \eqref{1D_Sequence_Condition}, from \eqref{1D_Fractal_Condition} we get $S = R^{\alpha /(2r+\alpha)}$, 
hence 
\begin{equation}
\frac{\| \sup_{n}|e^{it_n\Delta} g_R|   \|_{L^2(\mathcal H^\alpha)}}{\|g_R\|_2} \gtrsim R^{ \frac{1-\alpha}{2} + \frac{r(2\alpha - 1)}{2r+\alpha} },  
\end{equation}
which due to $S \geq R^{1/2}$ will be subject to the restriction $r \leq \alpha / 2$. 

To make this rigorous, as in previous sections, define the measure 
\begin{equation}
\mu_R(E) = \frac{|E \cap X_R|}{|X_R|} \,  \Big(\frac{R}{S^2}\Big)^\alpha.
\end{equation}
It suffices to prove that it is $\alpha$-dimensional. 
First observe that $Q \ll RS$ implies that the balls in $X_R$ are disjoint, 
so
$|X_R| = (Q/S^2)/S = R^\alpha/S^{1+\alpha}$, 
and hence 
\begin{equation}
\mu_R(E) = |E \cap X_R| \, S^{1-\alpha}.
\end{equation}
Now, for  any ball $B_r$ of radius $r > 0$, 
\begin{itemize}
	\item If $r \leq 1/S$, then $|B_r \cap X_R| \leq |B_r| = r$, so 
	\begin{equation}
	\mu_R(B_r) \leq r S^{1-\alpha} = r^\alpha \, (rS)^{1-\alpha} \leq r^\alpha. 
	\end{equation}
	
	\item If $1/S \leq r \leq R/Q$, then $|B_r \cap X_R| \leq 1/S$, so 
	\begin{equation}
	\mu_R(B_r) \leq \frac{1}{S} S^{1-\alpha} = \frac{1}{S^\alpha}  \leq r^\alpha. 
	\end{equation}

	\item If $R/Q \leq r \leq R/S^2$, then $B_r$ intersects at most $rQ/R$ intervals of $X_R$, so
	\begin{equation}
	\mu_R(B_r) \leq r \frac{Q}{R}\, \frac{1}{S} S^{1-\alpha} 
	= r \, \Big( \frac{S^2}{R}  \Big)^{1-\alpha}
	= r^\alpha \, \Big( r\frac{S^2}{R}  \Big)^{1-\alpha}
	\leq r^\alpha. 
	\end{equation}	
	
	\item If $R/S^2 \leq r \leq 1$, then $B_r$ may intersect all $X_R$, so 
	\begin{equation}
	\mu_R(B_r) \leq \frac{Q}{S^2} \, \frac{1}{S} \,  S^{1-\alpha} 
	= \Big( \frac{R }{S^2} \Big)^\alpha
	\leq r^\alpha. 
	\end{equation}	 
	
\end{itemize}
Hence, $\mu_R$ is $\alpha$-dimensional with constant $\langle \mu_R \rangle \simeq 1$,
so
\begin{equation} 
\frac{\| \sup_{n}|e^{it_n\Delta} g_R|   \|_{L^2(d\mu_R)}}{\|g_R\|_2} \gtrsim  S^{1/2} \, \mu(X_R)^{1/2} 
= S^{1/2} \,  \Big( \frac{R}{S^2 }\Big)^{\alpha/2}              
= R^{ \frac{1-\alpha}{2} + \frac{r(2\alpha - 1)}{2r+\alpha} },                                                          
\end{equation} 
where $S = R^{\alpha/(2r+\alpha)}$ from \eqref{1D_Sequence_Condition} and \eqref{1D_Fractal_Condition}. 
Hence, when $r \leq \alpha/2$, 
the $\alpha$-dimensional fractal maximal estimate cannot hold in $H^s$ for $s < \frac{1-\alpha}{2} + \frac{r(2\alpha - 1)}{2r+\alpha}$. 
Finally, when $r = \alpha/2$ the exponent is 1/4, so 
$s_c^{\operatorname{max}}(\alpha, \ell^{r,\infty}) \geq 1/4$ for all $r \geq \alpha /2$. 
$\qed$

\appendix

\section{A lemma to compute the measure of $X_R$}
\label{Appendix_LemmaAnChuPierce}
The following lemma
is an adaptation of Proposition 4.2 in \cite{ACP} by An, Chu and Pierce. 
We reproduce it in a slightly more general form than in \cite[Lemma 3.4]{EP_V}, 
and we refer to \cite[Lemma 3.4]{EP_V} and \cite[Proposition 5.2]{EP_V2} for a detailed proof. 
\begin{lem}\label{Lemma_AnChuPierce}
Let $Q \gg 1$. Define
\begin{equation}
X_Q = \bigcup_{\substack{ q \simeq Q \\ q \text{ odd }}} \bigcup_{\substack{ (p_1,p') \in [0,q)^n \\ \gcd(p_1,q)=1 }} B_1\Big( \frac{p_1}{q}, h_1(Q) \Big) \times B_{n-1} \Big( \frac{p'}{q}, h_2(Q) \Big)
\end{equation}
such that $h_1(Q) \simeq 1/Q^{\alpha}$ and $h_2(Q) \simeq 1/Q^\beta$ for some $\alpha, \beta \geq 1$. 
Then, 
\begin{equation}
|X_Q| \simeq \frac{Q^{n+1}h_1(Q) h_2(Q)^{n-1}}{1 + Q^{n+1}h_1(Q) h_2(Q)^{n-1}}.
\end{equation}
\end{lem}


\section*{Acknowledgements}

C.-H. Cho is supported by NRF-RS202300239774 and NRF-2022R1A4A1018904.

D. Eceizabarrena has received funding from the European Union’s Horizon Europe research and innovation programme under Marie Sklodowska Curie Actions with grant agreement 101104250 - TIDE, 
from the Simons Foundation Collaboration Grant on Wave Turbulence (Nahmod’s award ID 651469) and from the American Mathematical Society and the Simons Foundation under an AMS-Simons Travel Grant for the period 2022-2024.

\end{document}